\documentclass[11pt]{article}
\usepackage[a4paper,margin=1in]{geometry}
\usepackage[T1]{fontenc}
\usepackage{lmodern}
\usepackage{microtype}
\usepackage{amsmath,amssymb,amsthm,mathtools}
\usepackage[hidelinks,pdfusetitle]{hyperref}
\numberwithin{equation}{section}
\newtheorem{theorem}{Theorem}[section]
\newtheorem{lemma}[theorem]{Lemma}
\newtheorem{proposition}[theorem]{Proposition}
\newcommand{\e}{\mathrm e}
\newcommand{\one}{\mathbf1}
\newcommand{\Qcal}{\mathcal Q}
\newcommand{\D}{\mathcal D}
\newcommand{\I}{\mathcal I}
\newcommand{\II}{\mathcal J}
\newcommand{\III}{\mathcal K}
\newcommand{\Err}{\mathcal E}
\newcommand{\Sg}{\mathfrak S}
\newcommand{\abs}[1]{\left|#1\right|}
\title{Fixed-endpoint triple singular series along a dyadic family}
\author{Bahar U\u{g}urdo\u{g}an}
\date{}
\begin{document}
\maketitle
\begin{abstract}
Let $\mathfrak S$ denote the Hardy--Littlewood singular series. We determine
the logarithmic term in the sum over triples with endpoints $0$ and
$h=3\cdot2^k$. As $k\to\infty$, for each fixed $0<\varepsilon<1/9$ we have
\[
\sum_{1\le a<h}\mathfrak S(\{0,a,h\})
=\mathfrak S(\{0,h\})(h-\log h)
+O_\varepsilon((\log h)^{8/9+\varepsilon}).
\]
After removing the local factors at $2$ and $3$, the triple sum splits
into two additive convolutions. We estimate their centered parts using
reciprocal exponential sums and a second moment over products of divisors.
\end{abstract}

\section{Introduction and main result}\label{sec:intro}
For a set $\mathcal H$ of distinct integers put
\[
\mathfrak S(\mathcal H)=\prod_p
\left(1-\frac{\nu_p(\mathcal H)}p\right)
\left(1-\frac1p\right)^{-\#\mathcal H},
\]
where $\nu_p(\mathcal H)$ is the number of residue classes occupied by
$\mathcal H$ modulo $p$.  This singular series is the familiar local
factor in the Hardy--Littlewood prime-tuples problem \cite{HL}.
For even $h$, write $\Sigma(h)=\mathfrak S(\{0,h\})$ and
\[
 A_3(h)=\sum_{1\le a<h}\mathfrak S(\{0,a,h\}).
\]
Odlyzko, Rubinstein and Wolf \cite[Theorem~3]{ORW} obtained the leading
term $h\Sigma(h)$ in their work on jumping champions. Goldston and Ledoan
\cite[\S5]{GL} write $A_3(h)=h\Sigma(h)+R_3(h)$ and report the bound
$R_3(h)\ll_\varepsilon h^\varepsilon$. They also discuss numerical
evidence for a remainder of order at most $\Sigma(h)\log h$.
Feng and Wu \cite[Lemma~3.1]{FW} give a uniform leading-term estimate
when the interior shift ranges over a shorter interval.

Here we determine the coefficient of $\log h$ for endpoints whose
separation is three times a power of two. Related averages with
congruence restrictions or smooth weights are studied by
Kuperberg \cite{Kuperberg}. De la Bret\`eche \cite{deLaBreteche}
studies the modified triple singular series averaged over all three shifts.
The pair factor \(\Sigma(h)\) also appears in the endpoint-conditioned
sieve model for rough numbers in prime gaps~\cite[\S1.1]{Ugurdogan}.

Let $N=2^j$, $j\ge1$, $h=6N$, and $L=\log(2N)$.
The implied constants below are absolute unless a dependence is indicated.
Define
\[
\mathcal E(x)=x^{8/9}(\log(2x))^{5/2},\qquad x\ge1.
\]

\begin{theorem}\label{th:paper-main}
For $h=6\cdot2^j$, as $j\to\infty$,
\begin{equation}\label{eq:paper-main}
\sum_{1\le a<h}\mathfrak S(\{0,a,h\})
=\Sigma(h)(h-\log h)+O(\mathcal E(L)).
\end{equation}
In particular, for every fixed $0<\varepsilon<1/9$ the error is
$O_\varepsilon((\log h)^{8/9+\varepsilon})$.
\end{theorem}

On this family $\Sigma(h)=\Sigma(6)>0$, since the primes dividing $h$
are exactly $2$ and $3$. Thus $R_3(h)/(\Sigma(h)\log h)\to-1$.
In the notation of \cite[(2--7)]{ORW}, $A_{d,1}=\Sigma(2d)$ and
$A_{d,2}=-A_3(2d)$. Only even interior shifts contribute to $A_3(2d)$.
Theorem~\ref{th:convolution} bounds the centered convolution at each
of the two arguments $N$ and $3N$. The estimates use Suzuki's short
reciprocal sum \cite{Suzuki}, the Weil bound for complete Kloosterman
sums \cite{DFI}, and a bilinear moment inequality of the type used
by Bourgain and Garaev \cite{BG}.

\section{Arithmetic preliminaries and the dyadic difference}\label{sec:arithmetic}
Let $\Qcal$ be the set of positive squarefree integers coprime to $6$. Put
\[
\lambda(1)=1,\quad
\lambda(d)=\prod_{p\mid d}(p-3)^{-1}\ (d\in\Qcal),\quad
\lambda(d)=0\ (d\notin\Qcal),
\]
\[
f(n)=\sum_{d\mid n}\lambda(d),\qquad g(n)=f(n)-1,
\qquad C_b(M)=\sum_{1\le n<M}b(n)b(M-n).
\]
A divisor is called supported if it belongs to $\Qcal$.
With $t_p=1/[p(p-3)]$ for primes $p\ge5$, set
\[
K=\prod_{p\ge5}(1+t_p),\qquad P=\prod_{p\ge5}(1+2t_p),
\qquad \mu_\infty=P-2K+1.
\]
These products converge absolutely. For $M$ with no prime factor outside
$\{2,3\}$ write $Z(M)=C_g(M)-M\mu_\infty$.
Write $\e(t)=\exp(2\pi i t)$ and $\e_m(b)=\e(b/m)$, including
$\e_1(b)=1$ for integral $b$. For each positive integer $M$ define
\[
\delta(M)=C_g(2M)-2C_g(M).
\]

\begin{lemma}\label{lem:dyadic-inversion}
For $M=N$ and $M=3N$ we have
\begin{equation}\label{eq:dyadic-inversion}
 Z(M)=-\sum_{i\ge0}2^{-i-1}\delta(2^iM),
\end{equation}
with absolute convergence.  Moreover $Z(M)\ll\log^2(2M)$.
\end{lemma}
\begin{proof}
Since $g(2n)=g(n)$, separating even and odd indices gives
\[\delta(M)=\sum_{1\le n<2M}(-1)^{n+1}g(n)g(2M-n).\]
The sums $\sum_{d\le x}\lambda(d)\ll\log(2x)$ and
$\sum_{d>x}\lambda(d)/d\ll1/x$ for $x\ge1$ follow by dyadic summation
from \eqref{eq:weight-global} below. Expand $C_g(M)$ over $d,e>1$.
If $M$ has only $2,3$ as prime factors, a nonzero congruence count
requires $(d,e)=1$.  The Chinese remainder theorem gives
$M/(de)+O(1)$ solutions in the interior interval.  Summing the
errors over $d,e<M$ costs $O(\log^2(2M))$. Extending the main
term to all $d,e>1$ costs $O(1)$, since $\sum_d\lambda(d)/d=K$.
The Euler product
for the extended coefficient is $\mu_\infty$.
Iterating $\delta(M)=Z(2M)-2Z(M)$ and using the growth estimate gives
\eqref{eq:dyadic-inversion}.
\end{proof}

We use half-open dyadic ranges $d\sim D$ to mean $D<d\le2D$ and
$\D=\{2^i:0\le i\le j+2\}$.  Put
\[
s(x)=2\mathbf1_{[1/2,1)}(\{x\})-1,\qquad
U_{D,E}(M)=\sum_{\substack{d\sim D,e\sim E\\(d,e)=1}}
\lambda(d)\lambda(e)s(M\overline e/d).
\]
The multiplicative inverse $\overline e$ in this expression is taken
modulo $d$.  Set
$\mathcal B_{D,E}(M)=\tfrac12(U_{D,E}(M)+U_{E,D}(M))$.

\begin{lemma}\label{lem:inverse}
For $M=N$ or $M=3N$,
\begin{equation}\label{eq:delta-blocks}
\delta(M)=\sum_{D,E\in\D}U_{D,E}(M)
             =\sum_{D,E\in\D}\mathcal B_{D,E}(M).
\end{equation}
\end{lemma}
\begin{proof}
For coprime odd $d,e>1$ with $(M,de)=1$, let $V_{d,e}(2M)$
be the alternating sign sum over positive $n<2M$ with
$d\mid n$ and $e\mid 2M-n$. Let $B(t)$ count the positive solutions
of $dx+ey=t$. For $t=M,2M$ coprimality excludes solutions on the
axes, and Popoviciu's formula \cite{BR} gives
\[
 B(t)=1+\frac{t}{de}
       -\{t\overline e/d\}-\{t\overline d/e\}.
\]
At $t=2M$, $x$ and $y$ have the same parity, and the solutions with
both even are counted by $B(M)$. Thus
\[
V_{d,e}(2M)=B(2M)-2B(M)=\frac12\left\{
 s(M\overline e/d)+s(M\overline d/e)\right\}.
\]
Here $2\{x\}-\{2x\}=\mathbf1_{[1/2,1)}(\{x\})$.
In the divisor expansion of $g(n)g(2M-n)$ a common divisor of $d$ and
$e$ divides $2M$. Since $d,e$ are coprime to $6$, they must be
coprime to one another. Every pair with a solution has $d,e<2M\le6N$.
We may extend the symmetrized sum to the square $1<d,e\le8N$ since
each added pair has $B(2M)=B(M)=0$. Transposing the two variables
in this square gives
\eqref{eq:delta-blocks}.
\end{proof}

Let $\ell(n)=1/n$ and let $*$ denote Dirichlet convolution.
Then $\lambda=v*\ell$, where $v$ is multiplicative, with
\[
v(2)=-\frac12,\quad v(3)=-\frac13,\quad
v(p)=\frac3{p(p-3)},\quad v(p^2)=-\frac1{p(p-3)}\quad(p\ge5),
\]
all unlisted prime-power values being zero.
The absolute local factors imply
\begin{equation}\label{eq:weight-global}
\sum_{n\ge1}|v(n)|n^{1/4}<\infty,\quad
\sum_{d\sim D}\lambda(d)\ll1,\quad
\sum_{d\le x}\lambda(d)\ll\log(2x).
\end{equation}
The second bound follows by substituting $d=ak$ and summing $1/k$
over a dyadic interval. Summing over such intervals gives the third bound.

\section{Two divisor scales}\label{sec:scales}
In this section $M$ is either $N$ or $3N$.
For large $N$ put
\[
X=\exp\!\big(L^{2/3}\sqrt{\log L}\big).
\]
We partition the ordered dyadic pairs into the three disjoint sets
\begin{align*}
\I &=\{(D,E):E\ge D^{3/4}\},\\
\II&=\{(D,E):E<D^{3/4},\ DE\le N^{9/10}\},\\
\III&=\{(D,E):E<D^{3/4},\ DE>N^{9/10}\}.
\end{align*}

\begin{lemma}\label{lem:row}
For every $\epsilon>0$,
\[
 |U_{D,E}(M)|\ll_\epsilon
 \min\{1,(D^{1/2+\epsilon}/E)^{1/4}\}.
\]
\end{lemma}
\begin{proof}
Fix a supported modulus $q=d>1$ and a unit $b$ modulo $q$. On the units put $f(k)=s(b\bar k/q)$ and extend by zero elsewhere. Its period sum vanishes by $k\mapsto-k$, since $q$ is odd and $b$ is a unit. We claim that every interval sum of $f$ is $O_\epsilon(q^{1/2+\epsilon})$. Remove full periods, leaving a cyclic interval $I$. On residues $y\ne0$ set $s_q(y)=s(y/q)$ and set $s_q(0)=0$. For normalized Fourier coefficients the geometric sum formula gives
\[
 \sum_m|\widehat{\mathbf1_I}(m)|\ll\log(2q),\quad
 \widehat{s_q}(0)=0,\quad |\widehat{s_q}(n)|\ll 1/|n|
 \quad(0<|n|\le q/2).
\]
Indeed $s_q$ is a difference of two interval indicators. Expanding both factors in $\sum_{k\bmod q}^*\mathbf1_I(k)s_q(b\bar k)$ leaves complete Kloosterman sums $S(m,nb;q)$. The Weil bound \cite[Appendix, (6.3)]{DFI} gives
\[
 |S(m,nb;q)|\le \tau(q)\sqrt q\,(n,q)^{1/2},
\]
where $\tau$ counts divisors. For squarefree $q$ this also follows by multiplying the prime-modulus bounds, with the local factor $p-1$ if both coefficients vanish. Since
$\sum_{1\le n\le q}(n,q)^{1/2}/n\ll\log(2q)\sum_{r\mid q}r^{-1/2}$,
the standard elementary bound $\tau(q)\ll_\epsilon q^\epsilon$ proves the claim.

Expand $e=ak$ with $\lambda=v*\ell$, keeping $(a,q)=1$. In the $k$-sum use the claim with $b=M\bar a$, a unit since $(M,q)=1$. Partial summation and the trivial harmonic estimate bound the sum weighted by $1/k$ by
$O_\epsilon(\min\{1,q^{1/2+\epsilon}a/E\})$. Use $\min(1,x)\le x^{1/4}$ and \eqref{eq:weight-global}, then sum over $d$. The trivial bound follows separately from \eqref{eq:weight-global}.
\end{proof}

\begin{lemma}\label{lem:shift}
For \(R,F\ge4\), \(R\le A\le B\le2R\) and \(0\le\theta<1\), one has
\[
 \abs{\sum_{A<n+\theta\le B}\e\!\left(\frac{F}{n+\theta}\right)}
 \ll R\exp\!\left(-c\frac{(\log R)^3}{(\log F)^2}\right)+R^2/F
\]
with absolute constants.
\end{lemma}
\begin{proof}
By differencing two prefixes it suffices to take \(A=R\), \(B=R'\).
Put \(P=\lfloor R\rfloor\), \(X=\lfloor R-\theta\rfloor+1\) and
\(H=\lfloor R'-\theta\rfloor-\lfloor R-\theta\rfloor\). The sum begins
at \(X\) if \(H>0\), and \(0\le H\le\lceil R\rceil\le P+1\).
If \(P<3^{12}\), use the trivial bound \(R+1\). For \(F\le R\) the term
\(R^2/F\) suffices, and for \(F>R\) the exponential on the right is
at least \(\exp(-c\log R)\), bounded below on this range of \(R\).
Henceforth let \(P\ge3^{12}\).

First assume \(3P\le F^{2/3}\). We use Karatsuba's derivative estimate
\cite[Theorem~1]{Karatsuba} with the constants and normalization of
\cite[Lemma~1 and Remark~7]{Suzuki}. Its integer interval is \([X,X+P]\),
and its prefix length is \(P_1=\min(H,P)\). We may assume \(P_1>0\).
Dropping a possible last term costs at most one. On the whole interval,
\[
 P\le R<x+\theta\le R+1+P\le3P,\qquad
 D_m(x):=\frac{|f^{(m)}(x)|}{m!}=\frac{F}{(x+\theta)^{m+1}},
 \qquad f(x)=\frac{F}{x+\theta}.
\]
In particular \(f\) is real and smooth there. Write
\[
 a=\frac{\log F}{\log P},\quad \rho=\frac{2a}{3}>1,\quad
 K=\lfloor26a\rfloor,\quad \mathcal J=\mathbb Z\cap(2\rho,4\rho].
\]
Take \((c_0,c_1,c_2,c_3)=(1/39,25/26,23/24,1/4)\), and in the
original formulation take \(c_4=(c_1-c_2)/2=1/624\).
These satisfy \(0<c_0<1\), \(0<c_3\le c_2<c_2+c_4<c_1<1\).
We have \(K\ge2\) and
\[
 |\mathcal J|=\lfloor4\rho\rfloor-\lfloor2\rho\rfloor
 \ge2\rho-1\ge\rho\ge K/39.
\]
Each \(m\in\mathcal J\) satisfies \(2\le m\le4\rho\le K\).
For condition (A), \(a<(K+1)/26\) gives
\(D_{K+1}(x)\le P^{a-K-2}\le P^{-25(K+1)/26}\).
For condition (B), \(m>4a/3\) gives \(D_m(x)\le P^{a-m-1}\le P^{-m/4}\).
For the lower bound, \(F/(3P)\ge F^{1/3}\), \(3P\le P^{13/12}\)
and \(m\le8a/3\) give
\[
 D_m(x)\ge\frac{F}{3P}(3P)^{-m}
       \ge P^{a/3-13m/12}\ge P^{-23m/24}.
\]
All the derivative bounds hold on the same interval for each prefix.
The derivative estimate gives \(O(P^{1-\gamma/K^2})\) with
\(0<\gamma\le1\) and an implied constant depending only on the fixed
\(c_i\). For \(\log P\le C K^2\) with fixed \(C\), the trivial
\(O(P)\) bound already suffices since \(P^{\gamma/K^2}\le e^{C\gamma}\).
Either convention for including the ends of the integer sum
changes it by at most two terms, absorbed since \(P^{1-\gamma/K^2}\ge1\).
As \(K\le26\log F/\log P\) and \(\log P\ge\tfrac12\log R\),
this yields the asserted exponential term with \(c=\gamma/(8\cdot26^2)\).

If \(3P>F^{2/3}\), work on the real prefix \([R-\theta,R'-\theta]\).
The derivative of \(f\) is monotone there. For
\(\lambda_0=F/(8R^2)\) we have
\[
 0<\lambda_0< F/(4R^2)\le |f'(x)|\le F/R^2
 <3^{3/2}/\sqrt R<1/2<1-\lambda_0.
\]
For \(R'-R\ge1\) the Kusmin--Landau inequality \cite[Lemma~2]{Suzuki}
gives \(O(\lambda_0^{-1})=O(R^2/F)\).
A shorter prefix contains at most one integer and \(R^2/F>2\).
\end{proof}

\begin{lemma}\label{lem:middle}
For \((D,E)\in\II\), with absolute \(c>0\),
\begin{equation}\label{eq:decay}
 |U_{D,E}(M)|\ll \exp\!\left(-c\,\frac{(\log(2D))^3}{L^2}\right).
\end{equation}
\end{lemma}
\begin{proof}
Fix a supported \(e\sim E\). By \(e\bar e+d\bar d\equiv1\bmod de\),
\[
 s(M\bar e/d)=s\!\left(M/(de)-M\bar d/e\right).
\]
Expand \(d=ak\) using \(\lambda=v*\ell\), keeping \((a,e)=(k,e)=1\).
Reciprocity continues to hold for every such \(d>1\), including terms
of the expansion for which \((d,M)>1\).
The terms with \(a>D^{1/8}\) contribute \(O(D^{-1/32})\) per \(e\)
by \eqref{eq:weight-global}. For the other terms write \(k=en+r\) with \(1\le r<e\)
and \((r,e)=1\), and put
\[
 \theta=r/e,\quad R=D/(ae),\quad F=\nu M/(ae^2),\quad
 W=(\log(2D))^3/L^2 .
\]
Then \(R<n+\theta\le2R\), and the phase multiplied by the positive integer
\(\nu\) is \(F/(n+\theta)\) plus a constant on this class.
Since \(e<2D^{3/4}\), we have \(R>D^{1/8}/2\).
Bounded \(D\) follows from \eqref{eq:weight-global}. For sufficiently large
\(D\), \(\log R\gg\log(2D)\). We have
\[
 F/R=\nu M/(eD)\ge\tfrac12 N^{1/10}.
\]
Take \(H_0=\lceil\exp(W)\rceil\). Since \(D\le4N\), \(W\ll L\),
and \(\log F\ll L\) uniformly for \(1\le\nu\le H_0\).
Thus \(R,F\ge4\) for sufficiently large \(D\). Uniformly in the
residue \(r\) and in \(\nu\le H_0\),
\((\log R)^3/(\log F)^2\gg W\) and \(R/F\le2N^{-1/10}\).
Lemma~\ref{lem:shift} bounds each partial exponential sum divided
by \(R\) by \(O(\exp(-c_1W)+N^{-1/10})\).

Use the Erd\H os--Tur\'an inequality \cite[Ch.~2, Th.~2.5]{KN}
for the half-open interval \([1/2,1)\), on each initial segment of
the class. There are at most \(R+1\) points. The absolute sign sum divided
by \(R\) is
\[
 \ll H_0^{-1}+(1+\log H_0)(\exp(-c_1W)+N^{-1/10})
 \ll \exp(-c_2W).
\]
The last inequality uses \(W\ll L\) and a sufficiently small fixed
\(c_2>0\). Partial summation against \(1/k\le a/D\) costs at most a
constant times \(a/D\). The sum over the at most \(e\) unit classes
therefore costs \(O(1)\) in the normalization, since \(eRa/D=1\). Sum with \(|v(a)|\) and restore the tail, which is
\(O(\exp(-c_3W))\) since \(W\ll\log(2D)\). Finally sum over \(e\)
by \eqref{eq:weight-global}. All constants are independent of \(e,a,D,E,N\).
\end{proof}

Let \(\mu\) be the M\"obius function. For supported \(e>1\) set
\begin{align}
 b_e&=\sum_{(a,e)=1}v(a)
     =\frac13\prod_{\substack{p\ge5\\p\nmid e}}(1+2t_p),\nonumber\\
 m_e(u)&=\frac1e\sum_{b\mid e}\mu(b)s(u/b),\qquad
 J_e(D;M)=\int_D^{2D}m_e(M/x)\,\frac{dx}{x}.
 \label{eq:model}
\end{align}
The product for \(b_e\) follows from \eqref{eq:weight-global} and
\((1-1/2)(1-1/3)=1/3\).

\begin{lemma}\label{lem:model}
For each such \(e\) and every \(D\in\D\),
\[
 \sum_{\substack{d\sim D\\(d,e)=1}}\lambda(d)s(M\bar e/d)
 =b_e J_e(D;M)
       +O\!\left((e/D+M/D^2)^{1/4}\right).
\]
\end{lemma}
\begin{proof}
Use the same expansion \(d=ak\) and reciprocity as in
Lemma~\ref{lem:middle}, now for every \(a\) coprime to \(e\).
On each unit class \(k\equiv r\bmod e\), write
\(z_r\equiv M\overline{ar}\pmod e\).
Replace the sum over \(k\) on \(D/a<k\le2D/a\) by
\[
 \frac1e\int_{D/a}^{2D/a}s\!\left(\frac{M}{eat}-\frac{z_r}{e}\right)
                                                   \frac{dt}{t}.
\]
The sign has \(O(1+M/(eD))\) constant pieces because its argument is
monotone and its jumps are spaced by \(1/2\). On each piece the number
of progression points differs from the length divided by \(e\) by \(O(1)\)
for either endpoint convention. Partial summation against \(1/t\) gives
an error \(O(a/D)\) per piece, including possible endpoint terms. The sum
of the errors over the classes is \(O(a(e/D+M/D^2))\). It is also \(O(1)\)
by the dyadic harmonic bound and the integral bound \(\log 2\).

The \(z_r\) permute the units because \((Ma,e)=1\). For odd \(q\) one has
\(\sum_{r=0}^{q-1}s((w-r)/q)=s(w)\). The sum is periodic of period one
in \(w\), and on \(0\le w<1\) the number of positive signs increases
from \((q-1)/2\) to \((q+1)/2\) at \(w=1/2\), including that point.
Using \(\one_{(r,e)=1}=\sum_{b\mid(r,e)}\mu(b)\) gives exactly
\[
 \frac1e\sum_{r\bmod e}^{*}s((u-r)/e)=m_e(u).
\]
Thus the sum of the integrals, after \(x=at\), is \(J_e(D;M)\).
Sum with \(v(a)\), bounding the error by \(\min(1,y)\le y^{1/4}\) and
\eqref{eq:weight-global}. This also handles \(a>2D\), when the discrete sum is empty.
\end{proof}

\begin{lemma}\label{lem:large}
For each $M=N,3N$,
\[
\sum_{(D,E)\in\III}|U_{D,E}(M)|\ll1.
\]
\end{lemma}
\begin{proof}
If $(D,E)\in\III$, then $D^{7/4}>N^{9/10}$ and
$D>N^{18/35}$.  For $e\sim E$ the quantities $e/D$ and $M/D^2$
are $O(N^{-1/35})$.  Lemma~\ref{lem:model} replaces each
$U_{D,E}$ by $\sum_{e\sim E}\lambda(e)b_eJ_e(D;M)$ with
$O(N^{-1/140})$ error per block.  Because $s$ has mean zero and a
bounded primitive, a change of variable and integration by parts give
\[
\left|\int_D^{2D}s(M/(bx))\,\frac{dx}{x}\right|
\ll\min(1,bD/M).
\]
Summation over dyadic $D\le4N\le4M$ yields
\[
\sum_{D\in\D}|J_e(D;M)|\ll\tau(e)\log(2e)/e.
\]
The weighted sum of the right-hand side over $e>1$ is absolutely
bounded. Indeed $\tau(e)\le2\sqrt e$, and the dyadic
bound in \eqref{eq:weight-global} gives convergence.  There are
$O(L^2)$ dyadic pairs, so their total replacement error is $O(1)$.
\end{proof}

\begin{proposition}\label{prop:oriented-reduction}
For $M=N,3N$,
\begin{equation}\label{eq:reduce}
\delta(M)=\sum_{\substack{D,E\in\D\\D<X,\ E<D^{3/4}}}
U_{D,E}(M)+O(1).
\end{equation}
\end{proposition}
\begin{proof}
In $\I$, Lemma~\ref{lem:row} with $\epsilon=1/8$
gives $O(D^{-1/32})$ after summing in $E$ by a geometric series.
This is summable over dyadic $D$.  In $\III$ apply Lemma~\ref{lem:large}.
For $(D,E)\in\II$ with $D\ge X$, Lemma~\ref{lem:middle}
gives $\exp(-c(\log L)^{3/2})$ per block, and there are $O(L^2)$
blocks.  For $D<X$, $E<D^{3/4}$ we have $DE<X^2<N^{9/10}$
for sufficiently large $N$, so no further region remains.
\end{proof}

\section{Reciprocal energies of two primes}\label{sec:prime}
Let $P^+(r)$ denote the largest prime factor of $r$, with $P^+(1)=1$.
For $u=\log D$ and a dyadic prime interval $(P,2P]$ write
$t=\log(2P)$.  We shall use
\begin{equation}\label{eq:smooth}
\sum_{P^+(r)\le z}\lambda(r)\ll\log(2z),\qquad
\sum_{P^+(r)\le z}\lambda(r)\log r\ll(\log(2z))^2.
\end{equation}
The first estimate follows from Mertens' theorem applied to
$\prod_{5\le p\le z}(1+\lambda(p))$. For the second, multiply this
product by $\sum_{5\le p\le z}\lambda(p)\log p/(1+\lambda(p))\ll\log(2z)$.
Throughout the prime estimates below, $u$ is sufficiently large and
\begin{equation}\label{eq:ts}
1000\log u\le t\le0.76u.
\end{equation}

Put \(n_P=4P/t\) and \(k=k(P)=\lfloor0.96u/t\rfloor\). In what
follows we may take \(u\ge1000\). By the Brun--Titchmarsh inequality
with modulus one \cite[Theorem~2, (1.12)]{MV}, the prime count in
\((P,2P]\) is at most \(2P/\log P\le n_P\). Also
\(n_P=2e^t/t\le e^t\) and \(P\ge u^{1000}/2>k\ge1\).
Let \(J_k(P;d)\) count ordered \(2k\)-tuples of primes in \((P,2P]\),
all coprime to \(d\), for which the first and last \(k\) inverses have
the same sum modulo \(d\).

\begin{lemma}\label{lem:energy}
For each dyadic \(D\) there is a set \(\mathcal X_D\) of moduli \(d\sim D\)
with \(\sum_{d\in\mathcal X_D}\lambda(d)\ll D^{-1/50}\), outside
which, simultaneously for every dyadic \(P\) satisfying \eqref{eq:ts},
\[
 J_{k(P)}(P;d)\le (k(P)!+1)n_P^{\,k(P)}.
\]
The implied constant is absolute.
\end{lemma}
\begin{proof}
If two reciprocal sums are equal over the rationals, multiply their
difference by the product of the distinct primes in the tuple. Reducing
modulo each such prime shows that it divides its multiplicity difference.
Since all primes exceed \(k\), these differences vanish. The number of
these tuples is at most \(k!n_P^k\), including repetitions.

For every other tuple multiply the difference by the product of all
\(2k\) entries. The resulting nonzero integer \(T\) satisfies
\[
 |T|\le 2k(2P)^{2k}\le 2uD^{1.92}\le D^2.
\]
Every modulus counted for this tuple divides \(|T|\). Fix \(\epsilon=1/400\).
On the support of \(\lambda\),
\(\prod_{p\mid d}p/(p-3)\ll_\epsilon d^\epsilon\), since all sufficiently
large \(p\) satisfy \(p/(p-3)\le p^\epsilon\) and the small primes
occur at most once. The divisor bound now gives
\[
 \sum_{\substack{d\sim D\\d\mid T}}\lambda(d)
 \ll_\epsilon D^{-1+\epsilon}|T|^\epsilon
 \ll D^{-397/400}.
\]
The total weighted off-diagonal count is thus
\(O(D^{-397/400}n_P^{2k})\). The moduli where that count per modulus
exceeds \(n_P^k\) have total weight \(O(D^{-397/400}n_P^k)\) by
Markov's inequality. As \(n_P^k\le e^{kt}\le D^{384/400}\), this is
\(O(D^{-13/400})\). Take the union over at most
\(1+0.76u/\log2\) dyadic intervals. The asserted bound follows from
\((1+u)e^{-u/80}\ll1\). The only divisor estimate used a fixed \(\epsilon\).
\end{proof}

\begin{lemma}\label{lem:bilinear}
Suppose \(d\sim D\), \(d\notin\mathcal X_D\), and \((b_0,d)=1\).
Let \(P,Q\) satisfy \eqref{eq:ts}, with \(t_1=\log(2P)\),
\(t_2=\log(2Q)\), \(t_2\le t_1\) and \(t_1+t_2\le0.76u\).
Put \(v_0=t_1t_2/u\), \(J=\lceil e^{v_0/1000}\rceil\).
For every integer \(1\le\nu\le J\) and all \(|\alpha_p|,|\beta_q|\le1\),
\begin{equation}\label{eq:bilinear}
 \left|\sum_{p\sim P,\ q\sim Q}^{*}\alpha_p\beta_q
       \e_d(\nu b_0\bar p\bar q)\right|
       \le n_P n_Q e^{-v_0/50}.
\end{equation}
Here the variables are primes and the star imposes coprimality with \(d\).
\end{lemma}
\begin{proof}
Write \(k=k(P)\), \(l=k(Q)\), \(h_0=\nu b_0\), \(g_0=(h_0,d)=(\nu,d)\).
Let \(A_{P,k}(x)\) count the ordered \(k\)-tuples in the starred
\(p\)-sum with inverse sum \(x\bmod d\). Define \(A_{Q,l}(y)\) similarly.
Then \(\sum_x A_{P,k}(x)^2=J_k(P;d)\), and the analogous identity
holds for \(Q,l\).

We give the moment argument of \cite[\S4, (5)]{BG} allowing \(g_0>1\).
Call the sum in \eqref{eq:bilinear} \(S\) and write
\(B_p=\sum_q^*\beta_q\e_d(h_0\bar p\bar q)\).
H\"older gives \(|S|^l\le n_P^{l-1}\sum_p^*|B_p|^l\).
Choose \(|c_p|=1\) so that \(c_p B_p^l=|B_p|^l\). Expanding the powers,
grouping the \(q\)-tuples by \(y\), and using \(|\beta_q|\le1\) gives
\[
 |S|^l\le n_P^{l-1}\sum_y A_{Q,l}(y)|B(y)|,\qquad
 B(y)=\sum_p^*c_p\e_d(h_0\bar p y).
\]
Weighted H\"older with exponent \(k\), followed by Cauchy--Schwarz
and \(\sum_y A_{Q,l}(y)\le n_Q^l\), gives
\[
 |S|^{2kl}\le n_P^{2k(l-1)}n_Q^{2l(k-1)}
              J_l(Q;d)\sum_{y\bmod d}|B(y)|^{2k}.
\]
Let \(C_k(x)\) be the sum of \(\prod_i c_{p_i}\) over the tuples
counted by \(A_{P,k}(x)\). Orthogonality and Cauchy--Schwarz give
\[
 \sum_{y\bmod d}|B(y)|^{2k}
 =d\sum_{w\bmod d/g_0}
       \left|\sum_{\substack{x\bmod d\\x\equiv w\bmod d/g_0}} C_k(x)\right|^2
 \le d g_0\sum_{x\bmod d}A_{P,k}(x)^2.
\]
Indeed, each class modulo \(d/g_0\) contains exactly \(g_0\) classes
modulo \(d\). Thus, with constant one for every \(k,l\ge1\),
\begin{equation}\label{eq:moment}
 |S|^{2kl}\le dg_0\,n_P^{\,2kl-2k}n_Q^{\,2kl-2l}
                     J_k(P;d)J_l(Q;d).
\end{equation}

After Lemma~\ref{lem:energy} and division by \((n_Pn_Q)^{2kl}\), the
upper bound is \(\Delta=dg_0(k!+1)(l!+1)/(n_P^k n_Q^l)\).
The inequality \(\lfloor z\rfloor\ge z/2\) for \(z\ge1\) gives
\(kt_1\ge0.48u\). Since \(t_2\le0.38u\), we have
\(lt_2\ge0.96u-t_2\ge0.58u\). Hence \(kt_1+lt_2\ge1.06u\).
For \((m,t)=(k,t_1),(l,t_2)\),
\[
 m\log t\le u/1000,\qquad
 \log(m!+1)\le\log2+m\log m\le\log2+u/1000.
\]
Here \(m\le u/t\), \(t\ge1000\log u\), and both \(m,t\le u\).
Using \(\log n_P\ge t_1-\log t_1\) and its analogue for \(Q\),
\(d\le2e^u\), \(g_0\le J\le2e^{v_0/1000}\), and \(v_0\le u\), we get
\[
 \log\Delta
 \le 4\log2+u+v_0/1000+4u/1000-1.06u
 \le 4\log2-0.055u\le -u/25.
\]
Taking the \(2kl\)-th root and using \(kl\le u^2/(t_1t_2)=u/v_0\)
proves \eqref{eq:bilinear}.
\end{proof}

\begin{lemma}\label{lem:cut}
Under the hypotheses on \(d,b_0,P,Q\) of Lemma~\ref{lem:bilinear},
give each \(p\sim P\) any interval \(I(p)\) of \(q\)-values, possibly
empty. With \(v_0=t_1t_2/u\),
\[
 \left|\sum_{\substack{p\sim P,\ q\sim Q\\q\in I(p)}}^{*}
   \lambda(p)\lambda(q)s(b_0\bar p\bar q/d)\right|
       \ll \frac{e^{-v_0/2000}}{t_1t_2}
\]
with an absolute constant independent of the intervals.
\end{lemma}
\begin{proof}
On the integer lattice an interval indicator is a difference of two
prefix indicators. Consider one prefix \(q\le y_p\), clipping each
threshold to \([Q,2Q]\). Set \(x_j=Q+jQ/J\) for \(0\le j\le J\),
with \(J\) from Lemma~\ref{lem:bilinear}. Group the \(p\)'s by
\(j=\max\{i:x_i\le y_p\}\). The complete part of a group has the
common prefix \(q\le x_j\); at most \(J\) such parts are nonempty.
Since \(\lambda(p)\le2/P\), \(\lambda(q)\le2/Q\), their combined
contribution at each frequency \(1\le\nu\le J\) is bounded by
\(64J e^{-v_0/50}/(t_1t_2)\) by \eqref{eq:bilinear}.

For each \(p\) the incomplete part, if present, lies in one entire
bin \((x_j,x_{j+1}]\) of length \(Y=Q/J\).
As \(v_0\le0.76t_2\), \(\log J\le\log2+0.76t_2/1000\le\tfrac12\log Q\).
Thus \(Y\ge\sqrt Q>1\) and \(\log Y\ge t_2/4\).
By \cite[Theorem~2, (1.12)]{MV}, every entire bin contains at most
\(2Y/\log Y\le8Q/(Jt_2)\) primes. Multiplying by \(2/Q\) and
\(\sum_p^*\lambda(p)\le8/t_1\) bounds the combined incomplete parts by
\(128/(Jt_1t_2)\). These estimates also cover strict endpoints by
choosing the integer prefix thresholds before clipping.

The total weight is at most \(64/(t_1t_2)\). Apply Erd\H os--Tur\'an
\cite[Ch.~2, Th.~2.5]{KN} to the points \(b_0\bar p\bar q/d\bmod1\)
with \(q\in I(p)\), height \(J\), and test interval \([1/2,1)\).
If \(K_0\) is a common denominator, replace each point of weight \(w\)
by \(K_0w\) copies and divide the unnormalized inequality by \(K_0\).
The zero-mass case is immediate. The constant part in
\(s\) cancels since the interval has length \(1/2\). After multiplication
by \(t_1t_2\) the bound is
\[
 \ll J^{-1}+(1+\log J)(J e^{-v_0/50}+J^{-1})
 \ll (1+v_0/1000)e^{-v_0/1000}\ll e^{-v_0/2000},
\]
using \(e^{v_0/1000}\le J\le2e^{v_0/1000}\).
All constants here are absolute, including when \(0<v_0<1\).
\end{proof}

\section{A second moment at fixed divisor product}\label{sec:dispersion}
\begin{lemma}\label{lem:product}
Let \(R\ge2\), \(W\ge R^8\) and \(M=2^\alpha3^\beta\) with integers
\(\alpha,\beta\ge0\). For each supported \(1<e\le R\) fix an interval
\(J_e\subset(W,2W]\) of the form \((x_e,y_e]\), possibly empty. Then
\begin{equation}\label{eq:product}
 \left|\sum_{\substack{1<e\le R,\ e\in\Qcal,\ d\ge1\\de\in J_e,\ (d,e)=1}}
  \lambda(d)\lambda(e)s(M\bar e/d)\right|
 \ll \sqrt{\log(2R)}\,\bigl[\log(2+\log(2R))\bigr]^2
\end{equation}
with an absolute constant independent of \(M,R,W\) and of the intervals.
\end{lemma}
\begin{proof}
Write \(n=de\). On the support, \(\lambda(d)\lambda(e)=\lambda(n)\)
and the supported factorizations are precisely the divisors of a
supported \(n\). Here \(d>1\) since \(W\ge R^8\).
The reciprocity congruence
\(\bar e/d+\bar d/e\equiv1/n\pmod1\) gives, for every frequency
\(\nu\ge1\),
\begin{align*}
 S_\nu
 &=\sum_{\substack{1<e\le R,\ e\in\Qcal,\ de\in J_e\\(d,e)=1}}
       \lambda(d)\lambda(e)\e_d(\nu M\bar e)
   =\sum_{n\sim W}\lambda(n)\e(\nu M/n)F_\nu(n),\\
 F_\nu(n)
 &=\sum_{\substack{1<e\le R,\ e\mid n\\n\in J_e}}
                   \e_e(-\nu M\overline{n/e}).
\end{align*}
\(F_\nu\) is needed only on the support of \(\lambda\). By
\eqref{eq:weight-global}, the square of the absolute value of this sum is
\(\ll\sum_{n\sim W}\lambda(n)|F_\nu(n)|^2\).

For two divisors in the expansion put \(e_1=ba,\ e_2=bc\), where
\(b=(e_1,e_2)\), and let \(m=bac\). The factors \(b,a,c\) are
pairwise coprime and squarefree. The admissible products are \(n=mt\)
with \((t,m)=1\), weighted by \(\lambda(m)\lambda(t)\). If \(ac>1\),
the phase of the first divisor times the conjugate of the second is
\begin{equation}\label{eq:fold}
 \e_{ba}(-\nu M\overline{ct})\,
 \e_{bc}(\nu M\overline{at})=\e_{ac}(A\bar t),
 \quad A\equiv-\nu M\bar b\pmod a,\quad
       A\equiv \nu M\bar b\pmod c .
\end{equation}
Indeed its numerator over \(m\), before multiplication by \(\bar t\),
is \(\nu M(-c\bar c_{ba}+a\bar a_{bc})\), where the subscripts specify
the moduli. This numerator is zero modulo \(b\) and is congruent to
\(-\nu M\) modulo \(a\) and to \(\nu M\) modulo \(c\).
Since \((M,ac)=(b,ac)=1\), we have \((A,ac)=(\nu,ac)\).
This includes the cases \(a=1\) and \(c=1\) under the convention
that a congruence modulo one is vacuous. On the diagonal \(a=c=1\)
the phase is one.

Average this periodic function in \(t\), extended by zero off the
units modulo \(m\). Its mean equals
\begin{equation}\label{eq:ram}
 \frac{\varphi(b)}m c_{ac}(\nu),\qquad
 c_q(v)=\sum_{x\bmod q}^{*}\e_q(vx).
\end{equation}
Inversion permutes the units, and the Chinese remainder theorem
separates \(b\) from \(ac\). For squarefree \(q\) the local factors
of \(c_q(v)\) are \(p-1\) if \(p\mid v\) and \(-1\) otherwise,
as follows by summing the nonzero residues modulo each prime
\cite[Sec.~8.3]{Apostol}.
This also proves \(|c_q(v)|\le(q,v)\), including \(q=1\).

For an interval \((A_0,B_0]\subset(W/m,2W/m]\) we claim
\begin{equation}\label{eq:completion}
 \left|\sum_{\substack{A_0<t\le B_0\\(t,m)=1}}
             \lambda(t)\e_{ac}(A\bar t)\right|
 \ll \frac{(\nu,ac)}{ac}+\left(\frac{m^2}{W}\right)^{1/4}.
\end{equation}
Expand \(\lambda(t)\) by \eqref{eq:weight-global}, writing \(t=iy\)
with weight \(v(i)/y\).
Only \((i,m)=1\) occurs. The periodic function of \(y\) has the
same mean \(\alpha\) as in \eqref{eq:ram} since multiplication by
\(\bar i\) permutes the units. Its partial sums on any interval differ
from length times \(\alpha\) by \(O(m)\), by whole periods and two
end pieces. Partial summation with \(1/y\) gives the integral
\(\alpha\log(B_0/A_0)\) with error \(O(m^2i/W)\). Both the integral
and the harmonic sum are \(O(1)\) even when \(i>W/m\), so the error
is \(O(\min(1,m^2i/W))\). Summing with \(|v(i)|\), using
\(\min(1,z)\le z^{1/4}\), proves \eqref{eq:completion}, since
\(|\alpha|\le\varphi(b)(\nu,ac)/m\le(\nu,ac)/(ac)\) and
\(\log(B_0/A_0)\le\log2\). An empty
interval contributes zero. In the second moment the required
interval is \((J_{e_1}\cap J_{e_2})/m\), so the same estimate applies.

Since \(m\le R^2\) and \(W\ge R^8\), the total completion error is
\[
 \ll (R^4/W)^{1/4}
       \left(\sum_{i\le R}\lambda(i)\right)^3
   \ll R^{-1}(\log(2R))^3\ll1 .
\]
The three separate sums over \(b,a,c\le R\) dominate the sum of
\(\lambda(m)\). The contribution of the first term in
\eqref{eq:completion} is at most a constant times
\[
 \sum_{b\le R}\lambda(b)
 \left(\sum_{a\ge1}\frac{\lambda(a)(\nu,a)}a\right)^2
 \ll \log(2R)\,(\log(2\nu))^2 .
\]
For the last inequality, the Euler factors with \(p\nmid\nu\) have
bounded product, and the others are \(1+1/(p-3)\).
For \(p\ge5\), the ratio of \(1+1/(p-3)\) to
\((1-1/p)^{-1}\) is \(1+O(p^{-2})\).
Enlarging the product over \(p\mid\nu\) to all primes \(5\le p\le\nu\)
and applying Mertens' theorem \cite[Th.~3.4]{Kou} gives \(O(\log(2\nu))\),
also when \(\nu=1\). We have proved, uniformly for every \(\nu\ge1\),
\[
 \sum_{n\sim W}\lambda(n)|F_\nu(n)|^2
      \ll \log(2R)[\log(2\nu)]^2 .
\]
Consequently the original exponential sum \(S_\nu\) satisfies
\(|S_\nu|\ll\sqrt{\log(2R)}\log(2\nu)\).

Put \(T=\log(2R)\). The mass of the nonnegative weights is \(O(T)\), since
for each \(e\), the allowed \(d\)'s lie in \((W/e,2W/e]\), to which
\eqref{eq:weight-global} applies. Erd\H os--Tur\'an with
\(H=\lceil\sqrt T\rceil\) now bounds the sign sum in \eqref{eq:product} by
\[
 \ll \frac{T}{H+1}+
       \sqrt{T}\sum_{\nu=1}^H\frac{\log(2\nu)}{\nu}
 \ll \sqrt T\,[\log(2+T)]^2 .
\]
The half-open interval \([1/2,1)\) has length \(1/2\), so its constant
term cancels in \(s\). All bounds above are uniform for the prescribed
intervals, including their intersections in \eqref{eq:completion}.
\end{proof}

\section{The tail in the shorter divisor}\label{sec:tail}
For a dyadic row \(D\), put \(u=\log D\) and let
\(I_D=\bigcup_{E\in\D,\ E<D^{3/4}}(E,2E]\). It is a single
interval when nonempty.
\begin{lemma}\label{lem:tail}
Let \(M\in\{N,3N\}\). Uniformly for \(u\) sufficiently large and \(V\ge1\),
\[
 \left|\sum_{d\sim D}\lambda(d)
       \sum_{\substack{e\in I_D,\ \log e>V\\(e,d)=1}}
             \lambda(e)s(M\bar e/d)\right|
 \ll (\log(2u))^2+\frac{u\log(2u)}{V}.
\]
\end{lemma}
\begin{proof}
Put \(z_0=u^{1000}\). The supported \(e\le2D^{3/4}\) with at most
one prime factor above \(z_0\) have total weight at most
\[
 \prod_{5\le p\le z_0}(1+\lambda(p))
        \left(1+\sum_{z_0<p\le2D^{3/4}}\lambda(p)\right)
 \ll (\log(2u))^2 .
\]
Mertens' estimate gives this bound, and the sum of \(\lambda(d)\)
on the row is bounded. For other \(e\), write \(e=pqr\) uniquely
with \(p>q>z_0\) its two largest primes, and \(P^+(r)<q\), where
\(P^+(1)=1\). Take dyadic \(p\sim P_1,\ q\sim P_2\), and set
\(t_i=\log(2P_i)\). Thus \(t_2\le t_1\) and
\(t_1+t_2\le (3/4)u+\log8\le0.76u\) for large \(u\).

For fixed \(D,r,P_1,P_2\), the sum over \(d,p,q\), with the factor
\(\lambda(r)\) omitted, is
\begin{equation}\label{eq:prime-box}
 \ll (t_1t_2)^{-1}\exp(-c t_1t_2/u)
\end{equation}
in absolute value. Indeed, for fixed \(p,r\)
the conditions \(pqr\in I_D\), \(pqr>\exp V\),
\(P^+(r)<q<p\), \(q>z_0\) select an interval in \(q\).
These conditions are independent of \(d\). They define an
allowed interval in Lemma~\ref{lem:cut}, also when \(P_1=P_2\).
The lower size requirement \(t_i\ge1000\log u\) holds because both
primes exceed \(z_0\). For \((r,d)>1\) the contribution is zero. At every other
supported \(d\), \(b_0=M\bar r\) is a unit. Sum Lemma~\ref{lem:cut} over
the nonexceptional \(d\) with \eqref{eq:weight-global}. On the exceptional
set use the trivial prime bound \(O(1/(t_1t_2))\) and the total
\(\lambda\)-weight \(O(D^{-1/50})\) from Lemma~\ref{lem:energy}.
The same \(\mathcal X_D\) works for all \(P_1,P_2,r\).
Since \(v_0=t_1t_2/u\le0.76^2u/4\), one has
\(D^{-1/50}=e^{-u/50}\le e^{-v_0/2000}\).
Thus one may take \(c=1/2000\) in \eqref{eq:prime-box}.

A nonempty box has \(P^+(r)\le2P_2\) and
\(V<t_1+t_2+\log r\). In \eqref{eq:smooth} take \(z=2P_2\).
Hence the total \(\lambda(r)\) for the nonempty boxes with fixed
\(P_1,P_2\) is \(\ll ((t_1+t_2)t_2+t_2^2)/V\ll t_1t_2/V\).
This restriction is applied after \eqref{eq:prime-box}.
The \(t_i\) have spacing \(\log 2\). For each \(t_1\), the sum of
\(\exp(-c t_1t_2/u)\) over \(t_2\le t_1\) is
\(O(\min(t_1,u/t_1))\), by counting or a geometric series.
Summing to \(t_1\le0.76u\), splitting at \(\sqrt u\), costs
\(O(u\log(2u))\). This proves the lemma.
\end{proof}

\section{Bounds for the centered convolutions}\label{sec:convol}
\begin{theorem}\label{th:convolution}
For each $M=N,3N$,
\[
|C_g(2M)-2C_g(M)|+|Z(M)|\ll \mathcal E(L).
\]
\end{theorem}
\begin{proof}
Set $V=L^{4/9}$ and \(R=\exp V\). In \eqref{eq:reduce} all rows
\(D<R^{10}\) together cost \(O(V^2)\) by \eqref{eq:weight-global}, since
there are \(O((1+\log R)^2)\) ordered dyadic blocks. Restrict the
remaining rows to \(R^{10}\le D<X\). For \(e\le R\), partition the
products \(n=de\) into dyadic intervals \((W,2W]\). For each nonempty interval,
\(W>R^{10}/2\ge R^8\) and \(W<2XR\) for sufficiently large \(N\). For each fixed \(e\)
let \(E\) be its unique dyadic scale. The allowed \(D\)'s satisfy
\(R^{10}\le D<X,\ E<D^{3/4}\), hence form a consecutive set of
dyadic scales. Their \(d\)-intervals form a single half-open interval.
Multiplication by \(e\) and intersection with \((W,2W]\) gives
exactly an interval \(J_e\) as in Lemma~\ref{lem:product}.
In particular the product grouping precedes taking absolute values.
There are \(O(\log X)\) product intervals for large \(N\). This part is
\[
 \ll(\log X)\sqrt V\,[\log(2V)]^2
 \ll L^{8/9}(\log(2L))^{5/2}.
\]
For \(e>R\) sum Lemma~\ref{lem:tail} over the remaining \(D\).
Writing \(U=\log X=L^{2/3}\sqrt{\log L}\), the resulting bound is
\[
 \ll U[\log(2U)]^2+\frac{U^2\log(2U)}{V}
 \ll L^{2/3}(\log(2L))^{5/2}+L^{8/9}(\log(2L))^2.
\]
Since \(V^2=L^{8/9}\), we have proved \(\delta(M)\ll\Err(L)\).

Apply \eqref{eq:dyadic-inversion}. Replacing \(N\) by \(2^iN\)
preserves both endpoint families and replaces \(L\) by \(L+i\log2\).
The implied constant in the bound for \(\delta(2^iM)\) is unchanged, and
\(\Err(L+i\log2)\ll(1+i)^4\Err(L)\) for \(L\ge1\).
Absolute summation proves the assertion about \(Z(M)\).

\end{proof}

\section{Restoring the singular series}\label{sec:series}
To recover the linear terms in \(f=1+g\), put \(\sigma_{-1}(n)=\sum_{d\mid n}d^{-1}\).
For \(x\ge2\),
\begin{equation}\label{eq:sigma}
 \sum_{n\le x}\sigma_{-1}(n)=\zeta(2)x-\tfrac12\log x
                         +O((\log(2x))^{2/3}).
\end{equation}
Indeed, expand \(\sum_{t\le x}t^{-1}\lfloor x/t\rfloor\) and write
\(\psi(y)=\{y\}-1/2\). Extending \(x\sum_{t\le x}t^{-2}\) costs
\(O(1)\), and \(\sum_{t\le x}1/t=\log x+O(1)\), so it remains
to bound \(\sum_{t\le x}\psi(x/t)/t\).
For dyadic \(4\le A\le x\), put \(W_0=(\log A)^3/(\log(2x))^2\)
and \(H=\lceil\exp(W_0)\rceil\). On every prefix
\(A<t\le A'\le\min(2A,x)\), Suzuki's reciprocal estimate
\cite[Lemma 3]{Suzuki} gives, for \(1\le \nu\le H\),
\[
 \frac1A\left|\sum_{A<t\le A'}\e(\nu x/t)\right|
   \ll \exp(-cW_0)+\frac{A}{\nu x}.
\]
Here \(\log(\nu x)\ll\log(2x)\). Erd\H os--Tur\'an, followed by
integration using
\(\psi(y)=\int_0^1(\mathbf1_{[t,1)}(\{y\})-(1-t))\,dt\),
bounds the corresponding \(\psi\)-sum divided by \(A\) by
\(O(\exp(-c'W_0)+A/x)\). In the frequency sum the second term uses
\(\sum_{\nu\le H}\nu^{-2}\ll1\). Partial summation with \(1/t\)
preserves this last bound on the weighted row. The dyadic \(A/x\)
sum is bounded, and the exponential terms sum to
\(O((\log(2x))^{2/3})\) by integral comparison on the \(\log A\) grid.
The finitely many terms with \(t\le4\) cost \(O(1)\). This proves
\eqref{eq:sigma}.

Now \(f=v*\sigma_{-1}\), \(\sum_i v(i)=P/3\), and
\(\zeta(2)\sum_i v(i)/i=K\) by the local Euler factors of $v$. Convolving \eqref{eq:sigma} gives
\[
 \sum_{n\le x}g(n)=(K-1)x-\frac P6\log x+O((\log(2x))^{2/3}).
\]
Extending the linear coefficient costs \(O(x^{-1/4})\), since
\(x\sum_{i>x}|v(i)|/i\le x^{-1/4}\sum_i|v(i)|i^{1/4}\).
The series \(\sum |v(i)|\log(2i)\) converges, so the other endpoint
and coefficient extensions cost \(O(1)\). The estimate in
\eqref{eq:sigma} extends over \(1\le x<2\) by an absolute constant.
Both endpoints \(M=N,3N\) satisfy \(g(M)=0\), so the upper endpoint
may be omitted in the last sum. In
\(C_f(M)=M-1+2\sum_{0<n<M}g(n)+M\mu_\infty+Z(M)\),
the coefficient of \(M\) is \(1+2(K-1)+\mu_\infty=P\).
Theorem~\ref{th:convolution} gives
\begin{equation}\label{eq:cf}
 C_f(M)=PM-\frac P3\log M+O(\Err(L)).
\end{equation}

To evaluate the original triple sum, set
\[
 C_*=\prod_{p\ge5}\frac{p^2(p-3)}{(p-1)^3}.
\]
The product converges absolutely to a positive value.
An odd \(a\) gives zero at \(p=2\).
For \(a=2n\) the factors at \(2\) and \(3\) are \(4\) and
\((9/8)(1+\mathbf1_{3\mid n})\). For \(p\ge5\) the two endpoints
are distinct modulo \(p\). Relative to the generic factor in \(C_*\),
the correction is \(1+1/(p-3)\) if \(p\mid n(3N-n)\).
Such a prime cannot divide both factors. Hence, using \(f(3n)=f(n)\),
\[
 \sum_{1\le a<h}\Sg(\{0,a,h\})=\frac92 C_*\big(C_f(3N)+C_f(N)\big).
\]
The pair factors at \(2\) and \(3\) are \(2\) and \(3/2\), respectively.
For \(p\ge5\) one has
\[
 \frac{p^2(p-3)}{(p-1)^3}
 \left(1+\frac{2}{p(p-3)}\right)=\frac{p(p-2)}{(p-1)^2}.
\]
Thus \(\Sigma(h)=3C_*P\). By \eqref{eq:cf},
\[
 A_3(h)=\frac92 C_*P\left(4N-\frac{\log N+\log(3N)}3\right)
              +O(\Err(L)).
\]
Substitution of \(4N=2h/3\) and
\(\log N+\log(3N)=2\log h-\log12\) yields \eqref{eq:paper-main},
since the resulting constant \(\tfrac12\Sigma(6)\log12\) is absorbed
by \(O(\Err(L))\).

\end{document}